\documentclass{article}
\usepackage{amsmath}
\usepackage{amssymb}
\usepackage{amsfonts}
\usepackage{latexsym}
\usepackage{enumerate}

\newtheorem{llemma}{Lemma}[section]
\newtheorem{proposition}[llemma]{Proposition}

\newtheorem{example}[llemma]{Example}

\newtheorem{definition}[llemma]{Definition}
\newtheorem{remark}[llemma]{Remark}

\newtheorem{key}[llemma]{Keyword}

\usepackage[a4paper,left=2.5cm,right=2cm,top=3cm,bottom=2.5cm]{geometry}

\begin{document}

\title{\textbf{Relations on neutrosophic multi sets with properties}}
\author{Said Broumi$^{a}$, Irfan Deli$^{b}$, Florentin Smarandache$^{c}$ \\
~\\
$^{a}$ Faculty of Arts and Humanities, Hay El Baraka Ben M'sikCasablanca\\
B.P. 7951, Hassan II University Mohammedia-Casablanca , Morocco,\\
broumisaid78@gmail.com \\
~\\
$^{b}$ Muallim R{\i}fat Faculty of Education, \\
Kilis 7 Aral{\i}k University, 79000 Kilis, Turkey,\\
irfandeli@kilis.edu.tr \\
~\\
$^{c}$ Department of Mathematics, University of New Mexico,\\
705 Gurley Avenue, Gallup, NM 87301, USA \\
fsmarandache@gmail.com } \maketitle

\begin{abstract}

In this paper, we first give the cartesian product of two
neutrosophic multi sets(NMS). Then, we define relations on
neutrosophic multi sets to extend the intuitionistic fuzzy multi
relations to neutrosophic multi relations. The relations allows to
compose two neutrosophic sets. Also, various properties like
reflexivity, symmetry and transitivity are studied.
\end{abstract}

\textbf{\emph{AMS }} 03E72, 08A72

\begin{key}
Neutrosophic sets, neutrosophic multi sets, neutrosophic multi
relations, reflexivity, symmetry, transitivity.

\end{key}

\vspace{-3mm}
\noindent\hrulefill \\
\section{Introduction}

Recently, several theories have been proposed to deal with
uncertainty, imprecision and vagueness such as probability set theory,  fuzzy set theory\cite%
{zad-65}, intuitionistic fuzzy set theory \cite{ata-96}, rough set theory\cite%
{paw-82} etc. These theories are consistently being utilized as
efficient tools for dealing with diverse types of uncertainties and
imprecision embedded in a system. But, all these above theories
failed to deal with indeterminate and inconsistent information which
exist in beliefs system. In 1995, inspired from the sport games
(wining/tie/defeating), from votes (yes/ NA/ No), from decision
making (making a decision/ hesitating/not making) etc. and guided by
the fact that the law of excluded middle did not work any longer in
the modern logics, Smarandache\cite{sma-98} developed a new concept
called neutrosophic set (NS) which generalizes fuzzy sets and
intuitionistic fuzzy sets. NS can be described by membership degree,
indeterminate degree and non-membership degree. This theory and
their hybrid structures has proven useful in many different fields
such as control theory\cite{agg-10}, databases\cite{aro-10,aro-11},
medical diagnosis problem\cite{ans-11}, decision making problem
\cite{chi-13,kha-13}, physics\cite{rab-05}, topology \cite{lup-08},
etc. The works on neutrosophic set, in theories and
applications, have been progressing rapidly (e.g. \cite%
{ans-13,ash-02,bro-13a,wan-10}).

After Molodotsov\cite{mol-99} proposed the theory of soft set
combining fuzzy, intuitionistic fuzzy set models with other
mathematical models has attracted the attention of many researchers
(e.g. \cite{fen-14,muk-13,yan-13}. Also, Maji et al.\cite{maji-2013}
presented the concept of neutrosophic soft sets which is based on a
combination of the neutrosophic set and soft set models. Broumi and
Smarandache\cite{bro-13.,bro-13b} introduced the concept of the
intuitionistic neutrosophic soft set by combining the intuitionistic
neutrosophic sets and soft sets. The works on neutrosophic sets
combining soft sets, in theories and applications, have been
progressing rapidly (e.g. \cite{bro-13, bro-14,
bro-14a,del-14-1,del-13,maj-12}).

The notion of multisets was formulated first in \cite{yag-86} by
Yager as generalization of the concept of set theory and then the
multisets developed in \cite{cal-01} by Calude et al.. Several
authors from time to time made a number of generalization of set
theory. For example, Sebastian and Ramakrishnan\cite{seb-11a,seb-10}
introduced a new notion  called multi fuzzy sets, which is a
generalization of multiset. Since then, several
researcher\cite{mut-13,seb-11,syr-12,tho-13} discussed more
properties on multi fuzzy set. \cite{eje-14,shi-12} made an
extension of the concept of fuzzy multisets by an intuitionstic
fuzzy set, which called intuitionstic fuzzy multisets(IFMS). Since
then in the study on IFMS , a lot of excellent
results have been achieved by researchers \cite%
{das-13,raj-13,raj-13a,raj-13b,raj-14,raj-14a}. An element of a
multi fuzzy sets can occur more than once with possibly the same or
different membership values, whereas an element of intuitionistic
fuzzy multisets allows the repeated occurrences of membership and
non--membership values. The concepts
of FMS and IFMS fails to deal with indeterminatcy. In 2013 Smarandache \cite%
{sma-13} extended the classical neutrosophic logic to n-valued
refined neutrosophic logic, by refining each neutrosophic component
T, I, F into
respectively T$_{1}$, T$_{2}$, ..., T$_{m}$, and I$_{1}$, I$_{2}$, ..., I$%
_{p}$, and F$_{1}$, F$_{2}$, ..., F$_{r}$\textbf{.} Recently,  Ye et
al. \cite{ye-15}, Ye and Ye \cite{ye-15b} and Chatterjee et
al.\cite{cha-15} presented single valued neutrosophic multi sets in
detail. The concept of neutrosophic multi set (NMS)is a
generalisation of fuzzy multisets and intuitionistic fuzzy multi
sets.

The purpose of this paper is an attempt to extend the neutrosophic
relations to neutrosophic multi relations(NMR). This paper is
arranged in the following manner. In section 2, we present the basic
definitions and results of neutrosophic set theory and neutrosophic
multi(or refined) set theory that are useful for subsequent
discussions. In section 3, we study the concept of neutrosophic
multi relations and their operations. Finally, we conclude the
paper.
\section{Preliminary}
In this section, we present the basic definitions and results of
neutrosophic set theory \cite{sma-98,wan-10} and neutrosophic
multi(or refined) set theory \cite{del-14} that are useful for
subsequent discussions. See especially \cite
{aro-10,aro-11,ans-11,ans-13,ash-02,bro-13a,chi-13,del-14,del-14-1,kha-13,lup-08,rab-05}
for further details and background.

Smarandache\cite%
{sma-13} refine T , I, F to  $T_1$, $T_2$,..., $T_m$ and $I_1$,
$I_2$,..., $I_p$ and $F_1$, $F_2$,..., $F_r$  where all  $T_m$,
$I_p$ and  $F_r$ can be subset of [0,1]. In the following sections,
we considered only the case when T ,I and F are split into the same
number of subcomponents 1,2,...p, and $T_{A}^{j}$
$I_{A}^{j}$,$F_{A}^{j}$ are single valued neutrosophic number.
\begin{definition} \cite{sma-98}
Let U be a space of points (objects), with a generic element in U
denoted by u. A neutrosophic set(N-set) A in U is characterized by a
truth-membership function $T_A$, a indeterminacy-membership function
$I_A$ and a falsity-membership function $F_A$. $T_A(x)$, $I_A(x)$
and $F_A(x)$ are real standard or nonstandard subsets of
$]^-0,1^+[$.

It can be written as

$$A=\{<x,(T_A(x),I_A(x),F_A(x))>:x\in U, \,T_A(u),I_A(x),F_A(x)\subseteq [0,1]\}.$$

There is no restriction on the sum of $T_A(u)$; $I_A(u)$ and
$F_A(u)$, so $^-0\leq sup T_A(u) + sup I_A(u) + supF_A(u)\leq 3^+$.

Here, 1$^+$= 1+$\varepsilon$, where 1 is its standard part and
$\varepsilon$ its non-standard part. Similarly, $^-0$=
1+$\varepsilon$, where 0 is its standard part and $ \varepsilon $
its non-standard part.

\end{definition}
For application in real scientific and engineering areas,Wang et al.%
\cite{wan-10} proposed the concept of an SVNS, which is an instance
of neutrosophic set. In the following, we introduce the definition
of SVNS.
\begin{definition}
\cite{wan-10} Let U be a space of points (objects), with a generic
element in U denoted by u. An SVNS A inX is characterized by a
truth-membership function $T_A(x)$, a indeterminacy-membership
function $I_A(x)$ and a falsity-membership function $F_A(x)$, where
$T_A(x)$, $I_A(x)$, and $F_A(x)$ belongs to [0,1] for each point u
in U. Then, an SVNS A can be expressed as

\begin{equation*}
A=\{<u,(T_A(x),I_A(x),F_A(x))>:x\in E, \,T_A(x),I_A(x),F_A(x)\in
[0,1]\}.
\end{equation*}

There is no restriction on the sum of $T_A(x)$; $I_A(x)$ and $F_A(x)$, so $%
0\leq supT_A(x) +sup I_A(x) + supF_A(x)\leq 3$.
\end{definition}
\begin{definition}\cite{ye-15}
Let $E$ be a universe. A neutrosophic multi set (NMS or Nm-set) $A$
on $E$ can be defined as follows:
\begin{equation*}
\begin{array}{rr}
A & = \{<x,
(T_A^1(x),T_A^2(x),...,T_A^P(x)),(I_A^1(x),I_A^2(x),...,I_A^P(x)), \\
& (F_A^1(x),F_A^2(x),...,F_A^P(x))>: \,\, x \in E \}%
\end{array}%
\end{equation*}

where, $T_A^1(x),T_A^2(x),...,T_A^P(x): E \rightarrow [ 0, 1 ]$, $%
I_A^1(x),I_A^2(x),...,I_A^P(x): E \rightarrow [ 0, 1 ]$ and $%
F_A^1(x),F_A^2(x),...,F_A^P(x): E \rightarrow [ 0, 1 ]$ such that $0
\leq T_{A}^i (x) + I_A^i(x)+ F_A^i(x) \leq 3$($i=1,2,...,P$) and
$T_{A}^1(x)\leq
T_{A}^2(x)\leq...\leq T_{A}^P(x) $ for any $x\in E$. $%
(T_A^1(x),T_A^2(x),...,T_A^P(x))$, $(I_A^1(x),I_A^2(x),...,I_A^P(x))$ and $%
(F_A^1(x),F_A^2(x),...,F_A^P(x))$ is the truth-membership sequence,
indeterminacy-membership sequence and falsity-membership sequence of
the element $x$, respectively. Also, P is called the dimension of
NMS A. We arrange the truth-membership sequence in decreasing order
but the corresponding indeterminacy-membership and
falsity-membership sequence may not be in decreasing or increasing
order. \label{f} The Cardinality of the membership function Tc(x)
,the indterminacy function Ic(x) and non-membership Fc(x) is the the
lenght of an element x in a NMs A denoted as P(A), defined as

$$P(A)=\left\vert Tc(x)\right\vert =\left\vert Ic(x)\right\vert
=\left\vert Fc(x)\right\vert $$ if A,B,C are the NMS defined on E,
then their cardinality is $$P=Max[P(A),P(B),P(C)\}.$$
set of all Neutrosophic multi sets on E is denoted by NMS(E).
\end{definition}
\begin{definition}\cite{ye-15}
\label{x} Let $A,B \in NMS(E)$. Then,

\begin{enumerate}
\item $A$ is said to be NM subset of $B$ is denoted by $A\widetilde{\subseteq%
} B$ if $T_A^i(x)\leq T_B^i(x)$, $I_A^i(x)\geq I_B^i(x)$
,$F_A^i(x)\geq F_B^i(x)$, $\forall x\in E$.

\item $A$ is said to be neutrosophic equal of $B$ is denoted by $A= B$ if $%
T_A^i(x)= T_B^i(x)$, $I_A^i(x)= I_B^i(x)$, $F_A^i(x)= F_B^i(x)$,
$\forall x\in E$.
\item the complement of A denoted by $A^{\widetilde{c}}$ and is defined by
\begin{equation*}
\begin{array}{rr}
A^{\widetilde{c}} & = \{<x,(F_A^1(x),F_A^2(x),...,F_A^P(x))
,(1-I_A^1(x),1-I_A^2(x),...,1-I_A^P(x)), \\
& (T_A^1(x),T_A^2(x),...,T_A^P(x))>: \,\, x \in E \}%
\end{array}%
\end{equation*}
\end{enumerate}
\end{definition}

%
%


\begin{definition}\cite{ye-15}
\label{b} Let $A,B \in NMS(E)$. Then,

\begin{enumerate}
\item The union of $A$ and $B$ is denoted by $A\widetilde{\cup} B=C$ and
is defined by
\begin{equation*}
\begin{array}{rr}
C & = \{<x,
(T_C^1(x),T_C^2(x),...,T_C^P(x)),(I_C^1(x),I_C^2(x),...,I_C^P(x)), \\
& (F_C^1(x),F_C^2(x),...,F_C^P(x))>: \,\, x \in E \}%
\end{array}%
\end{equation*}
where $T_C^i=T_A^i(x)\vee T_B^i(x)$, $I_C^i=I_A^i(x)\wedge I_B^i(x)$ ,$%
F_C^i=F_A^i(x)\wedge F_B^i(x)$, $\forall x\in E$ and $i=1,2,...,P$.

\item The intersection of $A$ and $B$ is denoted by $A\widetilde{\cap} B=D$
and is defined by
\begin{equation*}
\begin{array}{rr}
D & = \{<x,
(T_D^1(x),T_D^2(x),...,T_D^P(x)),(I_D^1(x),I_D^2(x),...,I_D^P(x)), \\
& (F_D^1(x),F_D^2(x),...,F_D^P(x))>: \,\, x \in E \}%
\end{array}%
\end{equation*}
where $T_D^i=T_A^i(x)\wedge T_B^i(x)$, $I_D^i=I_A^i(x)\vee I_B^i(x)$ ,$%
F_D^i=F_A^i(x)\vee F_B^i(x)$, $\forall x\in E$ and $i=1,2,...,P$.

\item The addition of $A$ and $B$ is denoted by $A\widetilde{+} B={E_1}$ and
is defined by
\begin{equation*}
\begin{array}{rr}
{E_1} & = \{<x,
(T_{E_1}^1(x),T_{E_1}^2(x),...,T_{E_1}^P(x)),(I_{E_1}^1(x),I_{E_1}^2(x),...,I_{E_1}^P(x)),
\\
& (F_{E_1}^1(x),F_{E_1}^2(x),...,F_{E_1}^P(x))>: \,\, x \in E \}%
\end{array}%
\end{equation*}
where $T_{E_1}^i=T_A^i(x)+ T_B^i(x)-T_A^i(x). T_B^i(x)$, $%
I_{E_1}^i=I_A^i(x).I_B^i(x)$ ,$F_{E_1}^i=F_A^i(x). F_B^i(x)$,
$\forall x\in E $ and $i=1,2,...,P$.

\item The multiplication of $A$ and $B$ is denoted by $A\tilde{\times} B={E_2%
}$ and is defined by
\begin{equation*}
\begin{array}{rr}
{E_2} & = \{<x,
(T_{E_2}^1(x),T_{E_2}^2(x),...,T_{E_2}^P(x)),(I_{E_2}^1(x),I_{E_2}^2(x),...,I_{E_2}^P(x)),
\\
& (F_{E_2}^1(x),F_{E_2}^2(x),...,F_{E_2}^P(x))>: \,\, x \in E \}%
\end{array}%
\end{equation*}
where $T_{E_2}^i=T_A^i(x). T_B^i(x)$, $%
I_{E_2}^i=I_A^i(x)+I_B^i(x)-I_A^i(x).I_B^i(x)$ ,$F_{E_2}^i=F_A^i(x)+
F_B^i(x)-F_A^i(x). F_B^i(x)$, $\forall x\in E$ and $i=1,2,...,P$.
\end{enumerate}

\end{definition}
\section{Relations on Neutrosophic Multi Sets}

In this section, after given the cartesian product of two
neutrosophic  multi sets(NMS), we define relations on neutrosophic
multi sets and study their desired properties. The relation extend
the concept of intuitionistic multi relation \cite{raj-13b} to
neutrosophic multi relation. Some of it is quoted from
\cite{del-14,del-14-1,raj-13b,sma-98}.

\begin{definition}
Let $\emptyset \neq A,B\in NMS(E)$. Then, cartesian product of A and
B is a Nm-set in $E\times E$, denoted by $ A\times B$, defined as
\begin{equation*}
\begin{array}{rl}
A\times B= & \{<(x,y), (T_{A\times B}^1(x,y),T_{A\times B}^2(x,y),...,T_{A\times B}^n(x,y)), \\
& (I_{A\times B}^1(x,y),I_{A\times B}^2(x,y),...,I_{A\times B}^n(x,y)), \\
& (F_{A\times B}^1(x,y),F_{A\times B}^2(x,y),...,F_{A\times
B}^n(x,y))>: \,\, x,y \in E \}
\end{array}
\end{equation*}

where
\begin{equation*}
T_{A\times B}^{j},I_{A\times B}^{j},F_{A\times B}^{j}:E\rightarrow
\lbrack 0,1],
\end{equation*}%

\begin{equation*}
T_{A\times B}^{j}(x,y)=\min \left\{
T_{A}^{j}(x),T_{B}^{j}(y)\right\},
\end{equation*}

\begin{equation*}
I_{A\times B}^{j}(x,y)=\max \left\{
I_{A}^{j}(x),I_{B}^{j}(y)\right\}
\end{equation*}
and
\begin{equation*}
F_{A\times B}^{j}(x,y)=\max \left\{
F_{A}^{j}(x),F_{B}^{j}(y)\right\}
\end{equation*}

for all $x , y\in E$ and $j\in \{1,2,...,n\}$($n=max\{P(A),P(B)\})$.
\end{definition}

\begin{remark}
A cartesian product on A is a Nm-set in $E\times E$, denoted by
$A\times A$, defined as

\begin{equation*}
\begin{array}{rl}
A\times A= & \{<(x,y), (T_{A\times A}^1(x,y),T_{A\times A}^2(x,y),...,T_{A\times A}^n(x,y)), \\
& (I_{A\times A}^1(x,y),I_{A\times A}^2(x,y),...,I_{A\times A}^n(x,y)), \\
& (F_{A\times A}^1(x,y),F_{A\times A}^2(x,y),...,F_{A\times
A}^n(x,y))>: \,\, x,y \in E \}
\end{array}
\end{equation*}
where
\begin{equation*}
T_{A\times A}^{j},I_{A\times A}^{j},F_{A\times A}^{j}:E\times
E\rightarrow \lbrack 0,1],
\end{equation*}%

\begin{equation*}
T_{A\times A}^{j}(x,y)=\min \left\{
T_{A}^{j}(x),T_{A}^{j}(y)\right\},
\end{equation*}

\begin{equation*}
I_{A\times A}^{j}(x,y)=\max \left\{
I_{A}^{j}(x),I_{A}^{j}(y)\right\}
\end{equation*}
and
\begin{equation*}
F_{A\times A}^{j}(x,y)=\max \left\{
F_{A}^{j}(x),F_{A}^{j}(y)\right\}
\end{equation*}
 $j\in
\{1,2,...,n\}$($n=max\{P(A)\})$
\end{remark}


\begin{example}\label{1}
Let $E=\{x_1, x_2\}$ be a universal set and $A$ and $B$ be two
Nm-sets over E as;
$$
\begin{array}{cc}
  A= & \{<x_1,\{0.3,0.5,0.6\}, \{0.2,0.3,0.4\}
,\{0.4,0.5,0.9\}>, \\
    & <x_2,\{0.4,0.5,0.7\}, \{0.4,0.5,0.1\} ,\{0.6,0.2,0.7\}>\}
\end{array}
$$ and

$$
\begin{array}{cc}
  B= &\{<x_1,\{0.4,0.5,0.6\}, \{0.2,0.4,0.4\}
,\{0.3,0.8,0.4\}>,\\
    &<x_2,\{0.6,0.7,0.8\}, \{0.3,0.5,0.7\}
,\{0.1,0.7,0.6\}>\}
\end{array}
$$
 Then, the cartesian product of $A$ and $B$ is
obtained as follows

$$\begin{array}{rr}
A\times B =& \{<(x_1,x_1),\{0.3,0.5,0.6\}, \{0.2,0.4,0.4\}
,\{0.3,0.8,0.9\}>, \\& <(x_1,x_2),\{0.3,0.7,0.8\}, \{0.2,0.5,0.7\}
,\{0.1,0.7,0.9\}>,\\& <(x_2,x_1),\{0.4,0.5,0.6\}, \{0.2,0.5,0.4\}
,\{0.3,0.8,0.7\}>,\\& <(x_2,x_2),\{0.4,0.7,0.8\}, \{0.3,0.5,0.7\}
,\{0.1,0.7,0.7\}>\}
\end{array}
$$
\end{example}


\begin{definition}
Let $\emptyset \neq A,B\in NMS(E)$ and $j\in \{1,2,...,n\}$. Then, a
neutrosophic multi relation from $A$ to $B$ is a Nm-subset of
$A\times B$. In other words, a neutrosophic multi relation from $ A$
to $B$ is of the form $(R,C),\,(C\subseteq E\times E)$ where $
R(x,y)\subseteq A\times B$ $\forall (x,y)\in C$.
\end{definition}
\begin{example}\label{3}Let us consider the Example \ref{1}. Then, we define
a neutrosophic multi relation $R$ and $S$, from $A$ to $B$, as
follows
$$\begin{array}{rr}
R =& \{<(x_1,x_1),\{0.2,0.6,0.9\}, \{0.2,0.4,0.5\}
,\{0.3,0.8,0.9\}>,
\\& <(x_1,x_2),\{0.3,0.9,0.8\}, \{0.2,0.8,0.7\}
,\{0.1,0.8,0.9\}>,\\& <(x_2,x_1),\{0.1,0.9,0.6\}, \{0.2,0.5,0.4\}
,\{0.2,0.8,0.7\}>\}
\end{array}
$$

and

$$\begin{array}{rr}
S =& \{<(x_1,x_1),\{0.1,0.7,0.9\}, \{0.2,0.5,0.7\}
,\{0.1,0.9,0.9\}>,
\\& <(x_1,x_2),\{0.3,0.9,0.8\}, \{0.2,0.8,0.8\}
,\{0.1,0.8,0.9\}>,\\& <(x_2,x_1),\{0.1,0.9,0.7\}, \{0.2,0.9,0.4\}
,\{0.2,0.8,0.9\}>\}
\end{array}
$$
\end{example}

\begin{definition}
\label{b} Let $A,B\in NMS(E)$ and, R and S be two neutrosophic multi
relation from $A$ to $B$. Then, the operations $R\tilde{\cup}S$, $R\tilde{%
\cap}S$, $R\tilde{+}S$ and ${R\tilde{\times}S}$ are defined as
follows;

\begin{enumerate}
\item
\begin{equation*}
\begin{array}{rl}
R\widetilde{\cup} S = & \{<(x,y), (T_{R\widetilde{\cup} S}^1(x,y),T_{R%
\widetilde{\cup} S}^2(x,y),...,T_{R\widetilde{\cup} S}^n(x,y)), \\
& (I_{R\widetilde{\cup} S}^1(x,y),I_{R\widetilde{\cup} S}^2(x,y),...,I_{R%
\widetilde{\cup} S}^n(x,y)), \\
& (F_{R\widetilde{\cup} S}^1(x,y),F_{R\widetilde{\cup} S}^2(x,y),...,F_{R%
\widetilde{\cup} S}^n(x,y))>: \,\, x,y \in E \}%
\end{array}%
\end{equation*}
where

\begin{equation*}
T_{R\widetilde{\cup} S}^j(x,y)=T_R^j(x)\vee T_S^j(y),
\end{equation*}
\begin{equation*}
I_{R\widetilde{\cup} S}^j(x,y)=I_R^j(x)\wedge I_S^j(y),
\end{equation*}
\begin{equation*}
F_{R\widetilde{\cup} S}^j(x,y)=F_R^j(x)\wedge F_S^j(y)
\end{equation*}
$\forall x,y \in E$ and $j=1,2,...,n$.

\item
\begin{equation*}
\begin{array}{rl}
R\tilde{\cap} S & = \{<(x,y), (T_{R\widetilde{\cap} S}^1(x,y),T_{R\widetilde{%
\cap} S}^2(x,y),...,T_{R\widetilde{\cap} S}^n(x,y)), \\
& (I_{R\widetilde{\cap} S}^1(x,y),I_{R\widetilde{\cap} S}^2(x,y),...,I_{R%
\widetilde{\cap} S}^n(x,y)), \\
& (F_{R\widetilde{\cap} S}^1(x,y),F_{R\widetilde{\cap} S}^2(x,y),...,F_{R%
\widetilde{\cap} S}^n(x,y))>: \,\, x,y \in E \}%
\end{array}%
\end{equation*}
where

\begin{equation*}
T_{R\widetilde{\cap} S}^j(x,y)=T_R^j(x)\wedge T_S^j(y),
\end{equation*}
\begin{equation*}
I_{R\widetilde{\cap} S}^j(x,y)=I_R^j(x)\vee I_S^j(y),
\end{equation*}
\begin{equation*}
F_{R\widetilde{\cap} S}^j(x,y)=F_R^j(x)\vee F_S^j(y)
\end{equation*}
$\forall x,y \in E$ and $j=1,2,...,n$.

\item
\begin{equation*}
\begin{array}{rr}
{R\widetilde{+} S}= & \{<(x,y), (T_{{R\widetilde{+} S}}^1(x,y),T_{{R%
\widetilde{+} S}}^2(x,y),...,T_{{R\widetilde{+} S}}^n(x,y)), \\
& (I_{{R\widetilde{+} S}}^1(x,y),I_{{R\widetilde{+} S}}^2(x,y),...,I_{{R%
\widetilde{+} S}}^n(x,y)), \\
& (F_{{R\widetilde{+} S}}^1(x,y),F_{{R\widetilde{+} S}}^2(x,y),...,F_{{R%
\widetilde{+} S}}^n(x,y))>: \,\, x,y \in E \}%
\end{array}%
\end{equation*}
where

\begin{equation*}
T_{{R\widetilde{+} S}}^j(x,y)=T_R^j(x)+ T_S^j(y)-T_R^j(x). T_S^j(y),
\end{equation*}
\begin{equation*}
I_{{R\widetilde{+} S}}^j(x,y)=I_R^j(x).I_S^j(y),
\end{equation*}
\begin{equation*}
F_{{R\widetilde{+} S}}^j(x,y)=F_R^j(x). F_S^j(y)
\end{equation*}
$\forall x,y\in E$ and $j=1,2,...,n$.

\item
\begin{equation*}
\begin{array}{rr}
{R\tilde{\times}S}= & \{<(x,y), (T_{{R\tilde{\times}S}}^1(x,y),T_{R\tilde{%
\times}S}^2(x,y),...,T_{{R\tilde{\times}S}}^n(x,y)), \\
& (I_{{R\tilde{\times}S}}^1(x,y), I_{{R\tilde{\times}S}}^2(x,y),...,I_{{R%
\tilde{\times}S}}^n(x,y)), \\
& (F_{{R\tilde{\times}S}}^1(x,y),F_{{R\tilde{\times}S}}^2(x,y),...,F_{{R%
\tilde{\times}S}}^n(x,y))>: \,\, x,y \in E \}%
\end{array}%
\end{equation*}
where

\begin{equation*}
T_{{R\tilde{\times}S}}^j(x,y)=T_R^j(x). T_S^j(y),
\end{equation*}
\begin{equation*}
I_{{R\tilde{\times}S}}^j(x,y)=I_R^j(x)+I_S^j(y)-I_R^j(x).I_S^j(y),
\end{equation*}
\begin{equation*}
F_{{R\tilde{\times}S}}^j(x,y)=F_R^j(x)+ F_S^j(y)-F_R^j(x). F_S^j(y)
\end{equation*}
$\forall x,y\in E$ and $j=1,2,...,n$.
\end{enumerate}

Here $\vee$, $\wedge$, $+$, $.$, $-$ denotes maximum, minimum,
addition, multiplication, subtraction of real numbers respectively.
\end{definition}

%
%
\begin{example}\label{4} Let us consider the Example \ref{3}. Then,
$$\begin{array}{rr}
R\tilde{\cup} S =& \{<(x_1,x_1),\{0.2,0.6,0.9\}, \{0.2,0.4,0.5\}
,\{0.3,0.8,0.9\}>,
\\& <(x_1,x_2),\{0.3,0.9,0.8\}, \{0.2,0.8,0.7\}
,\{0.1,0.8,0.9\}>,\\& <(x_2,x_1),\{0.1,0.9,0.6\}, \{0.2,0.5,0.4\}
,\{0.2,0.8,0.7\}>\}
\end{array}
$$

and

$$\begin{array}{rr}
R\tilde{\cap} S =& \{<(x_1,x_1),\{0.1,0.7,0.9\}, \{0.2,0.5,0.7\}
,\{0.1,0.9,0.9\}>,
\\& <(x_1,x_2),\{0.3,0.9,0.8\}, \{0.2,0.8,0.8\}
,\{0.1,0.8,0.9\}>,\\& <(x_2,x_1),\{0.1,0.9,0.7\}, \{0.2,0.9,0.4\}
,\{0.2,0.8,0.9\}>\}
\end{array}
$$
Similarly, ${R\widetilde{+} S}$ and ${R\widetilde{\times} S}$ can be
computed.
\end{example}

Assume that $\emptyset \neq A,B,C\in NMS(E)$. Two neutrosophic multi
relations under a suitable composition, could too yield a new
neutrosophic multi relation with a useful significance. Composition
of relations is important for applications, because of the reason
that if a relation on A and B is known and if a relation on B and C
is known then the relation on A and C could be computed and defined
as follows;

\begin{definition}
Let R(A$\rightarrow $ \ B) and S (B$\rightarrow $\ C) be two
neutrosophic multi relations. The composition S $\circ $R is a
neutrosophic multi relation from A to C, defined by
\begin{equation*}
\begin{array}{rr}
{S\circ R}= & \{<(x,z),(T_{{S\circ R}}^{1}(x,z),T_{S\circ R}^{2}(x,z),...,T_{%
{S\circ R}}^{n}(x,z)), \\
& (I_{{S\circ R}}^{1}(x,z),I_{{S\circ R}}^{2}(x,z),...,I_{{S\circ R}%
}^{n}(x,z)), \\
& (F_{{S\circ R}}^{1}(x,z),F_{{S\circ R}}^{2}(x,z),...,F_{{S\circ R}%
}^{n}(x,z))>:\,\,x,z\in E\}%
\end{array}%
\end{equation*}%
where
\begin{equation*}
T_{S\circ R}^{j}(x,z)=\underset{y}{\vee }\left\{
T_{R}^{j}(x,y)\wedge T_{S}^{j}(y,z)\right\}
\end{equation*}

\begin{equation*}
I_{S\circ R}^{j}(x,z)=\underset{y}{\wedge }\left\{
I_{R}^{j}(x,y)\vee I_{S}^{j}(y,z)\right\}
\end{equation*}
and
\begin{equation*}
F_{S\circ R}^{j}(x,z)=\underset{y}{\wedge }\left\{
F_{R}^{j}(x,y)\vee F_{S}^{j}(y,z)\right\}
\end{equation*}

for every (x, z) $E\times E$, for every $y \in E$ and $j=1,2,...,n$.
\end{definition}


\begin{definition}
A neutrosophic multi relation R on A is said to be;

\begin{enumerate}
\item Reflexive if $T_{R}^{j}(x,x)=1$, $I_{R}^{j}(x,x)=0$ and $%
F_{R}^{j}(x,x)=0$ for all $x\in E$,

\item Symmetric if $T_{R}^{j}(x,y)=T_{R}^{j}(y,x)$, $%
I_{R}^{j}(x,y)=I_{R}^{j}(y,x)$ and $F_{R}^{j}(x,y)=F_{R}^{j}(y,x)$ for all $%
x,y \in E$,

\item Transitive if $R\circ R\subseteq R,$

\item neutrosophic multi equivalence relation if the relation R satisfies
reflexive, symmetric and transitive.
\end{enumerate}
\end{definition}


\begin{definition}
The transitive closure of a neutrosophic multi relation R on
$E\times E$
is $\overset{\symbol{94}}{R}=R\tilde{\cup}R^{2}\tilde{\cup}R^{3}\tilde{\cup}%
...$
\end{definition}

\begin{definition}
If R is a neutrosophic multi relation from A to B then $R^{-1}$ is
the inverse neutrosophic multi relation R from B to A, defined as
follows:

\begin{equation*}
R^{-1} =\left\{ \left\langle
(y,x),T_{R^{-1}}^{j}(x,y)),I_{R^{-1}}^{j}(x,y),F_{R^{-1}}^{j}(x,y)\right%
\rangle :(x,y)\in E\times E\right\}
\end{equation*}
where

$T_{R^{-1}}^{j}(x,y)=T_{R}^{j}(y,x),$ $I_{R^{-1}}^{j}(x,y)=I_{R}^{j}(y,x)$, $%
F_{R^{-1}}^{j}(x,y)=F_{R}^{j}(y,x)$ and $j=1,2,...,n$.
\end{definition}


\begin{proposition}
If R and S are two neutrosophic multi relation from A to B and B to
C, respectively. Then,

\begin{enumerate}
\item $(R^{-1})^{-1}=R$

\item $(S\circ R)^{-1}=R^{-1}\circ S^{-1}$
\end{enumerate}
\end{proposition}

\textbf{Proof}

\begin{enumerate}
\item Since $R^{-1}$ is a neutrosophic multi relation from B to A, we have

$T_{R^{-1}}^{j}(x,y)=T_{R}^{j}(y,x)$,
$I_{R^{-1}}^{j}(x,y)=I_{R}^{j}(y,x)$ and
$F_{R^{-1}}^{j}(x,y)=F_{R}^{j}(y,x)$

Then,

\begin{equation*}
T_{(R^{-1})^{-1}}^{j}(x,y)=T_{R^{-1}}^{j}(y,x)=T_{R}^{j}(x,y),
\end{equation*}%
\begin{equation*}
I_{(R^{-1})^{-1}}^{j}(x,y)=I_{R^{-1}}^{j}(y,x)=I_{R}^{j}(x,y)
\end{equation*}%
and
\begin{equation*}
F_{(R^{-1})^{-1}}^{j}(x,y)=F_{R^{-1}}^{j}(y,x)=F_{R}^{j}(x,y)
\end{equation*}%
therefore $(R^{-1})^{-1}=R$.

\item If the composition $S\circ R$ is a neutrosophic multi relation from
A to C, then the compostion $R^{-1}\circ S^{-1}$ is a neutrosophic
multi relation from C to A. Then,

\begin{equation*}
\begin{array}{rl}
T_{(S\circ R)^{-1}}^{j}(z,x) & =T_{(S\circ R)}^{j}(x,z) \\
& =\underset{y}{\vee }\left\{ T_{R}^{j}(x,y)\wedge T_{S}^{j}(y,z)\right\} \\
& =\underset{y}{\vee }\left\{ T_{R^{-1}}^{j}(y,x)\wedge
T_{S^{-1}}^{j}(z,y)\right\} \\
& =\underset{y}{\vee }\left\{ T_{S^{-1}}^{j}(z,y)\wedge
T_{R^{-1}}^{j}(y,x)\right\} \\
& =T_{R^{-1}\circ S^{-1}}^{j}(z,x)%
\end{array}%
,
\end{equation*}

\begin{equation*}
\begin{array}{rl}
I_{(S\circ R)^{-1}}^{j}(z,x) & =I_{(S\circ R)}^{j}(x,z) \\
& =\underset{y}{\wedge }\left\{ I_{R}^{j}(x,y)\vee I_{S}^{j}(y,z)\right\} \\
& =\underset{y}{\wedge }\left\{ I_{R^{-1}}^{j}(y,x)\vee
I_{S^{-1}}^{j}(z,y)\right\} \\
& =\underset{y}{\wedge }\left\{ I_{S^{-1}}^{j}(z,y)\vee
I_{R^{-1}}^{j}(y,x)\right\} \\
& =I_{R^{-1}\circ S^{-1}}^{j}(z,x)%
\end{array}%
\end{equation*}%
and
\begin{equation*}
\begin{array}{rl}
F_{(S\circ R)^{-1}}^{j}(z,x) & =F_{(S\circ R)}^{j}(x,z) \\
& =\underset{y}{\wedge }\left\{ F_{R}^{j}(x,y)\vee F_{S}^{j}(y,z)\right\} \\
& =\underset{y}{\wedge }\left\{ F_{R^{-1}}^{j}(y,x)\vee
F_{S^{-1}}^{j}(z,y)\right\} \\
& =\underset{y}{\wedge }\left\{ F_{S^{-1}}^{j}(z,y)\vee
F_{R^{-1}}^{j}(y,x)\right\} \\
& =F_{R^{-1}\circ S^{-1}}^{j}(z,x)%
\end{array}%
\end{equation*}
\end{enumerate}

Finally; proof is valid.

\begin{proposition}
If R is symmetric, then R$^{-1}$is also symmetric.
\end{proposition}

\textbf{Proof:} Assume that R is Symmetric then we have
\begin{equation*}
T_{R}^{j}(x,y)=T_{R}^{j}(y,x),
\end{equation*}%
\begin{equation*}
I_{R}^{j}(x,y)=I_{R}^{j}(y,x)
\end{equation*}%
and
\begin{equation*}
F_{R}^{j}(x,y)=F_{R}^{j}(y,x)
\end{equation*}

Also if R$^{-1}$ is an inverse relation, then we have
\begin{equation*}
T_{R^{-1}}^{j}(x,y)=T_{R}^{j}(y,x),
\end{equation*}
\begin{equation*}
I_{R^{-1}}^{j}(x,y)=I_{R}^{j}(y,x)
\end{equation*}
and
\begin{equation*}
F_{R^{-1}}^{j}(x,y)=F_{R}^{j}(y,x)
\end{equation*}
for all $x,y \in E$

To prove R$^{-1}$ is symmetric, it is enough to prove

\begin{equation*}
T_{R^{-1}}^{j}(x,y)=T_{R^{-1}}^{j}(y,x),
\end{equation*}
\begin{equation*}
I_{R^{-1}}^{j}(x,y)=I_{R^{-1}}^{j}(y,x)
\end{equation*}
and
\begin{equation*}
F_{R^{-1}}^{j}(x,y)=F_{R^{-1}}^{j}(y,x)
\end{equation*}
for all $x,y \in E$

Therefore;
\begin{equation*}
T_{R^{-1}}^{j}(x,y)=T_{R}^{j}(y,x)=T_{R}^{j}(x,y)=T_{R^{-1}}^{j}(y,x);
\end{equation*}
\begin{equation*}
I_{R^{-1}}^{j}(x,y)=I_{R}^{j}(y,x)=I_{R}^{j}(x,y)=I_{R^{-1}}^{j}(y,x)
\end{equation*}
and
\begin{equation*}
F_{R^{-1}}^{j}(x,y)=F_{R}^{j}(y,x)=F_{R}^{j}(x,y)=F_{R^{-1}}^{j}(y,x)
\end{equation*}
Finally; proof is valid.

\begin{proposition}
If R is symmetric ,if and only if $R= R^{-1}$.
\end{proposition}

\textbf{Proof:} Let R be symmetric, then
\begin{equation*}
T_{R}^{j}(x,y)=T_{R}^{j}(y,x);
\end{equation*}
\begin{equation*}
I_{R}^{j}(x,y)=I_{R}^{j}(y,x)
\end{equation*}
and
\begin{equation*}
F_{R}^{j}(x,y)=F_{R}^{j}(y,x)
\end{equation*}
and

R$^{-1}$ is an inverse relation, then
\begin{equation*}
T_{R^{-1}}^{j}(x,y)=T_{R}^{j}(y,x);
\end{equation*}
\begin{equation*}
I_{R^{-1}}^{j}(x,y)=I_{R}^{j}(y,x)
\end{equation*}
and
\begin{equation*}
F_{R^{-1}}^{j}(x,y)=F_{R}^{j}(y,x)
\end{equation*}
for all $x,y \in E$

Therefore; $T_{R^{-1}}^{j}(x,y)=T_{R}^{j}(y,x)=T_{R}^{j}(x,y).$

Similarly
\begin{equation*}
I_{R^{-1}}^{j}(x,y)=I_{R}^{j}(y,x)=I_{R}^{j}(x,y)
\end{equation*}
and
\begin{equation*}
F_{R^{-1}}^{j}(x,y)=F_{R}^{j}(y,x)=F_{R}^{j}(x,y)
\end{equation*}
for all $x,y \in E.$

Hence ${R=R}^{-1}$

Conversely, assume that ${R=R}^{-1}$ then, we have
\begin{equation*}
T_{R}^{j}(x,y)=T_{R^{-1}}^{j}(x,y)=T_{R}^{j}(y,x).
\end{equation*}
Similarly
\begin{equation*}
I_{R}^{j}(x,y)=I_{R^{-1}}^{j}(x,y)=I_{R}^{j}(y,x)
\end{equation*}
and
\begin{equation*}
F_{R}^{j}(x,y)=F_{R^{-1}}^{j}(x,y)=F_{R}^{j}(y,x).
\end{equation*}
Hence R is symmetric.

\begin{proposition}
If R and S are symmetric neutrosophic multi relations, then

\begin{enumerate}
\item $R\tilde{\cup} S$,

\item $R\tilde{\cap} S$,

\item $R \tilde{+}S$

\item ${R\tilde{\times}S}$
\end{enumerate}

are also symmetric.
\end{proposition}

\textbf{Proof:} R is symmetric, then we have;
\begin{equation*}
T_{R}^{j}(x,y)=T_{R}^{j}(y,x),
\end{equation*}
\begin{equation*}
I_{R}^{j}(x,y)=I_{R}^{j}(y,x)
\end{equation*}
and
\begin{equation*}
F_{R}^{j}(x,y)=F_{R}^{j}(y,x)
\end{equation*}
similarly S is symmetric, then we have
\begin{equation*}
T_{S}^{j}(x,y)=T_{S}^{j}(y,x),
\end{equation*}
\begin{equation*}
I_{S}^{j}(x,y)=I_{S}^{j}(y,x)
\end{equation*}
and
\begin{equation*}
F_{S}^{j}(x,y)=F_{S}^{j}(y,x)
\end{equation*}
Therefore,

\begin{enumerate}
\item
\begin{equation*}
\begin{array}{rl}
T_{R\widetilde{\cup} S}^{j}(x,y) & =\max \left\{
T_{R}^{j}(x,y),T_{S}^{j}(x,y)\right\} \\
& =\max \left\{ T_{R}^{j}(y,x),T_{S}^{j}(y,x)\right\} \\
& =T_{R\widetilde{\cup} S}^{j}(y,x)%
\end{array}%
,
\end{equation*}

\begin{equation*}
\begin{array}{rl}
I_{R\widetilde{\cup} S}^{j}(x,y) & =min\left\{
I_{R}^{j}(x,y),I_{S}^{j}(x,y)\right\} \\
& =\min \left\{ I_{R}^{j}(y,x),I_{S}^{j}(y,x)\right\} \\
& =I_{R\widetilde{\cup} S}^{j}(y,x),%
\end{array}%
\end{equation*}
and
\begin{equation*}
\begin{array}{rl}
F_{R\widetilde{\cup} S}^{j}(x,y) & =\min \left\{
F_{R}^{j}(x,y),F_{S}^{j}(x,y)\right\} \\
& =\min \left\{ F_{R}^{j}(y,x),F_{S}^{j}(y,x)\right\} \\
& =F_{R\widetilde{\cup} S}^{j}(y,x)%
\end{array}%
\end{equation*}
therefore, $R\widetilde{\cup} S$ is symmetric.

\item
\begin{equation*}
\begin{array}{rl}
T_{R\widetilde{\cap} S}^{j}(x,y) & =\min \left\{
T_{R}^{j}(x,y),T_{S}^{j}(x,y)\right\} \\
& =\min \left\{ T_{R}^{j}(y,x),T_{S}^{j}(y,x)\right\} \\
& =T_{R\widetilde{\cap} S}^{j}(y,x),%
\end{array}%
\end{equation*}

\begin{equation*}
\begin{array}{rl}
I_{R\widetilde{\cap} S}^{j}(x,y) & =max\left\{
I_{R}^{j}(x,y),I_{S}^{j}(x,y)\right\} \\
& =\max \left\{ I_{R}^{j}(y,x),I_{S}^{j}(y,x)\right\} \\
& =I_{R\widetilde{\cap} S}^{j}(y,x),%
\end{array}%
\end{equation*}
and
\begin{equation*}
\begin{array}{rl}
F_{R\widetilde{\cap} S}^{j}(x,y) & =\max \left\{
F_{R}^{j}(x,y),F_{S}^{j}(x,y)\right\} \\
& =\max \left\{ F_{R}^{j}(y,x),F_{S}^{j}(y,x)\right\} \\
& =F_{R\widetilde{\cap} S}^{j}(y,x)%
\end{array}%
\end{equation*}
therefore; $R\widetilde{\cap} S$ is symmetric.

\item
\begin{equation*}
\begin{array}{rl}
T_{R\tilde{+}S}^{j}(x,y) &
=T_{R}^{j}(x,y)+T_{S}^{j}(x,y)-T_{R}^{j}(x,y)T_{S}^{j}(x,y) \\
& =T_{R}^{j}(y,x)+T_{S}^{j}(y,x)-T_{R}^{j}(y,x)T_{S}^{j}(y,x) \\
& =T_{R\tilde{+}S}^{j}(y,x)%
\end{array}%
\end{equation*}

\begin{equation*}
\begin{array}{rl}
I_{R\tilde{+}S}^{j}(x,y) & =I_{R}^{j}(x,y)I_{S}^{j}(x,y) \\
& =I_{R}^{j}(y,x)I_{S}^{j}(y,x) \\
& =I_{R\tilde{+}S}^{j}(y,x)%
\end{array}%
\end{equation*}
and
\begin{equation*}
\begin{array}{rl}
F_{R\tilde{+}S}^{j}(x,y) & =F_{R}^{j}(x,y)F_{S}^{j}(x,y) \\
& =F_{R}^{j}(y,x)F_{S}^{j}(y,x) \\
& =F_{R\tilde{+}S}^{j}(y,x)%
\end{array}%
\end{equation*}

therefore, $R\tilde{+}S$ is also symmetric

\item
\begin{equation*}
\begin{array}{rl}
T_{R\tilde{\times} S}^{j}(x,y) & =T_{R}^{j}(x,y)T_{S}^{j}(x,y) \\
& =T_{R}^{j}(y,x)T_{S}^{j}(y,x) \\
& =T_{R \tilde{\times}t S}^{j}(y,x)%
\end{array}%
\end{equation*}

\begin{equation*}
\begin{array}{rl}
I_{R\tilde{\times} S}^{j}(x,y) &
=I_{R}^{j}(x,y)+I_{S}^{j}(x,y)-I_{R}^{j}(x,y)I_{S}^{j}(x,y) \\
& =I_{R}^{j}(y,x)+I_{S}^{j}(y,x)-I_{R}^{j}(y,x)I_{S}^{j}(y,x) \\
& =I_{R\tilde{\times} S}^{j}(y,x)%
\end{array}%
\end{equation*}

\begin{equation*}
\begin{array}{rl}
F_{R\tilde{\times} S}^{j}(x,y) &
=F_{R}^{j}(x,y)+F_{S}^{j}(x,y)-F_{R}^{j}(x,y)F_{S}^{j}(x,y) \\
& =F_{R}^{j}(y,x)+F_{S}^{j}(y,x)-F_{R}^{j}(y,x)F_{S}^{j}(y,x) \\
& =F_{R\tilde{\times}S}^{j}(y,x)%
\end{array}%
\end{equation*}

hence, $R\tilde{\times}S$ is also symmetric.
\end{enumerate}


\begin{remark}
R$\circ $S in general is not symmetric, as
\begin{equation*}
\begin{array}{rl}
T_{(R\circ S)}^{j}(x,z) & =\underset{y}{\vee }\left\{
T_{S}^{j}(x,y)\wedge
T_{R}^{j}(y,z)\right\} \\
& =\underset{y}{\vee }\left\{ T_{S}^{j}(y,x)\wedge T_{R}^{j}(z,y)\right\} \\
& \neq T_{(R\circ S)}^{j}(z,x)%
\end{array}%
\end{equation*}

\begin{equation*}
\begin{array}{rl}
I_{(R\circ S)}^{j}(x,z) & =\underset{y}{\wedge }\left\{
I_{S}^{j}(x,y)\vee
I_{R}^{j}(y,z)\right\} \\
& =\underset{y}{\wedge }\left\{ I_{S}^{j}(y,x)\vee I_{R}^{j}(z,y)\right\} \\
& \neq I_{(R\circ S)}^{j}(z,x)%
\end{array}%
\end{equation*}

\begin{equation*}
\begin{array}{rl}
F_{(R\circ S)}^{j}(x,z) & =\underset{y}{\wedge }\left\{
F_{S}^{j}(x,y)\vee
F_{R}^{j}(y,z)\right\} \\
& =\underset{y}{\wedge }\left\{ F_{S}^{j}(y,x)\vee F_{R}^{j}(z,y)\right\} \\
& \neq F_{(R\circ S)}^{j}(z,x)%
\end{array}%
\end{equation*}

but $R\circ S$ is symmetric, if $R\circ S=S\circ R$, for R and S are
symmetric relations.

\begin{equation*}
\begin{array}{rl}
T_{(R\circ S)}^{j}(x,z) & =\underset{y}{\vee }\left\{
T_{S}^{j}(x,y)\wedge
T_{R}^{j}(y,z)\right\} \\
& =\underset{y}{\vee }\left\{ T_{S}^{j}(y,x)\wedge T_{R}^{j}(z,y)\right\} \\
& =\underset{y}{\vee }\left\{ T_{R}^{j}(y,x)\wedge T_{R}^{j}(z,y)\right\} \\
& T_{(R\circ S)}^{j}(z,x)%
\end{array}%
\end{equation*}

\begin{equation*}
\begin{array}{rl}
I_{(R\circ S)}^{j}(x,z) & =\underset{y}{\wedge }\left\{
I_{S}^{j}(x,y)\vee
I_{R}^{j}(y,z)\right\} \\
& =\underset{y}{\wedge }\left\{ I_{S}^{j}(y,x)\vee I_{R}^{j}(z,y)\right\} \\
& =\underset{y}{\wedge }\left\{ I_{R}^{j}(y,x)\vee I_{R}^{j}(z,y)\right\} \\
& I_{(R\circ S)}^{j}(z,x)%
\end{array}%
\end{equation*}
and
\begin{equation*}
\begin{array}{rl}
F_{(R\circ S)}^{j}(x,z) & =\underset{y}{\wedge }\left\{
F_{S}^{j}(x,y)\vee
F_{R}^{j}(y,z)\right\} \\
& =\underset{y}{\wedge }\left\{ F_{S}^{j}(y,x)\vee F_{R}^{j}(z,y)\right\} \\
& =\underset{y}{\wedge }\left\{ F_{R}^{j}(y,x)\vee F_{R}^{j}(z,y)\right\} \\
& F_{(R\circ S)}^{j}(z,x)%
\end{array}%
\end{equation*}
for every $(x,z)\in E\times E$ and for $y\in E.$
\end{remark}


\begin{proposition}
If R is transitive relation, then $R^{-1}$ is also transitive.
\end{proposition}

\textbf{Proof :} R is transitive relation, if $R\circ R\subseteq R,$
hence if $R^{-1} \circ R^{-1}\subseteq R^{-1},$ then $R^{-1}$ is
transitive.

Consider;
\begin{equation*}
\begin{array}{rl}
T_{R^{-1}}^{j}(x,y) & =T_{R}^{j}(y,x)\geq T_{R\circ R}^{j}(y,x) \\
& = \underset{z}{\vee }\left\{ T_{R}^{j}(y,z)\wedge
T_{R}^{j}(z,x)\right\}
\\
& =\underset{z}{\vee }\left\{ T_{R^{-1}}^{j}(x,z)\wedge
T_{R^{-1}}^{j}(z,y)\right\} \\
& =T_{R^{-1}\circ R^{-1}}^{j}(x,y)%
\end{array}%
\end{equation*}

\begin{equation*}
\begin{array}{rl}
I_{R^{-1}}^{j}(x,y) & =I_{R}^{j}(y,x)\leq I_{R\circ R}^{j}(y,x) \\
& =\underset{z}{ \wedge }\left\{ I_{R}^{j}(y,z)\vee
I_{R}^{j}(z,x)\right\}
\\
& =\underset{z }{\wedge }\left\{ I_{R^{-1}}^{j}(x,z)\vee
I_{R^{-1}}^{j}(z,y)\right\} \\
& =I_{R^{-1}\circ R^{-1}}^{j}(x,y)%
\end{array}%
\end{equation*}
and
\begin{equation*}
\begin{array}{rl}
F_{R^{-1}}^{j}(x,y) & =F_{R}^{j}(y,x)\leq F_{R\circ R}^{j}(y,x) \\
& =\underset{z}{ \wedge }\left\{ F_{R}^{j}(y,z)\vee
F_{R}^{j}(z,x)\right\}
\\
& =\underset{z }{\wedge }\left\{ F_{R^{-1}}^{j}(x,z)\vee
F_{R^{-1}}^{j}(z,y)\right\} \\
& =F_{R^{-1}\circ R^{-1}}^{j}(x,y)%
\end{array}%
\end{equation*}
hence, proof is valid.

\begin{proposition}
If R is transitive relation, then $R \cap S$ is also transitive.
\end{proposition}

\textbf{Proof: } As R and S are transitive relations, $R \circ
R\subseteq R$ and $S\circ S\subseteq S$.

Also
\begin{equation*}
\begin{array}{cc}
T_{R\widetilde{\cap} S}^{j}(x,y)\geq T_{(R\widetilde{\cap} S)\circ (R%
\widetilde{\cap} S)}^{j}(x,y) &  \\
I_{R\widetilde{\cap} S}^{j}(x,y)\leq I_{(R\widetilde{\cap} S)\circ (R%
\widetilde{\cap} S)}^{j}(x,y) &  \\
F_{R\widetilde{\cap} S}^{j}(x,y)\leq F_{(R\widetilde{\cap} S)\circ (R%
\widetilde{\cap} S)}^{j}(x,y) &
\end{array}%
\end{equation*}
implies $R\widetilde{\cap} S)\circ (R\widetilde{\cap} S)\subseteq R
\cap S$, hence $R \cap S$ is transitive.

\begin{proposition}
If R and S are transitive relations, then

\begin{enumerate}
\item $R\tilde{\cup}S$,

\item $R \tilde{+}S$

\item ${R\tilde{\times}S}$
\end{enumerate}

are not transitive.
\end{proposition}

\textbf{Proof: }

\begin{enumerate}
\item As
\begin{equation*}
\begin{array}{cc}
T_{R\widetilde{\cup} S}^{j}(x,y)= \max \left\{
T_{R}^{j}(x,y),T_{S}^{j}(x,y)\right\} &  \\
I_{R\widetilde{\cup} S}^{j}(x,y)=\min \left\{
I_{R}^{j}(x,y),I_{S}^{j}(x,y)\right\} &  \\
F_{R\widetilde{\cup} S}^{j}(x,y)=\min \left\{
F_{R}^{j}(x,y),F_{S}^{j}(x,y)\right\} &
\end{array}%
\end{equation*}
and
\begin{equation*}
\begin{array}{cc}
T_{(R\widetilde{\cup} S)\circ (R\widetilde{\cup} S)}^{j}(x,y)\geq T_{R%
\widetilde{\cup} S}^{j}(x,y) &  \\
I_{(R\widetilde{\cup} S)\circ (R\widetilde{\cup} S)}^{j}(x,y)\leq I_{R%
\widetilde{\cup} S}^{j}(x,y) &  \\
F_{(R\widetilde{\cup} S)\circ (R\widetilde{\cup} S)}^{j}(x,y)\leq F_{R%
\widetilde{\cup} S}^{j}(x,y) &
\end{array}%
\end{equation*}

\item As
\begin{equation*}
\begin{array}{cc}
T_{R\tilde{+}%
S}^{j}(x,y)=T_{R}^{j}(x,y)+T_{S}^{j}(x,y)-T_{R}^{j}(x,y)T_{S}^{j}(x,y) &  \\
I_{R\tilde{+}S}^{j}(x,y)=I_{R}^{j}(x,y)I_{S}^{j}(x,y) &  \\
F_{R\tilde{+}S}^{j}(x,y)=F_{R}^{j}(x,y)F_{S}^{j}(x,y) &
\end{array}%
\end{equation*}
and

\begin{equation*}
\begin{array}{cc}
T_{(R\tilde{+}S)\circ (R\tilde{+}S)}^{j}(x,y)\geq
T_{R\tilde{+}S}^{j}(x,y) &
\\
I_{(R\tilde{+}S)\circ (R\tilde{+}S)}^{j}(x,y)\leq
I_{R\tilde{+}S}^{j}(x,y) &
\\
F_{(R\tilde{+}S)\circ (R\tilde{+}S)}^{j}(x,y)\leq
F_{R\tilde{+}S}^{j}(x,y) &
\end{array}%
\end{equation*}

\item As
\begin{equation*}
\begin{array}{cc}
T_{R\tilde{\times} S}^{j}(x,y)=T_{R}^{j}(x,y)T_{S}^{j}(x,y) &  \\
I_{R\tilde{\times}
S}^{j}(x,y)=I_{R}^{j}(x,y)+I_{S}^{j}(x,y)-I_{R}^{j}(x,y)I_{S}^{j}(x,y) &  \\
F_{R\tilde{\times}
S}^{j}(x,y)=F_{R}^{j}(x,y)+F_{S}^{j}(x,y)-F_{R}^{j}(x,y)F_{S}^{j}(x,y)
&
\end{array}%
\end{equation*}
and
\begin{equation*}
\begin{array}{cc}
T_{(R\tilde{\times} S)\circ (R\tilde{\times} S)}^{j}(x,y)\geq T_{R\tilde{%
\times} S}^{j}(x,y) &  \\
I_{(R\tilde{\times} S)\circ (R\tilde{\times} S)}^{j}(x,y)\leq I_{R\tilde{%
\times} S}^{j}(x,y) &  \\
F_{(R\tilde{\times} S)\circ (R\tilde{\times} S)}^{j}(x,y)\leq F_{R\tilde{%
\times} S}^{j}(x,y) &
\end{array}%
\end{equation*}
\end{enumerate}

Hence $R\tilde{\cup}S$, $R \tilde{+}S$ and ${R\tilde{\times}S}$ are
not
transitive. 

\begin{proposition}
If R is transitive relation, then $R^{2}$ is also transitive.
\end{proposition}

\textbf{Proof:} R is transitive relation, if $R \circ R\subseteq R$,
therefore if $R^{2} \circ R^{-2}\subseteq R^{2}$, then $R^{2}$ is
transitive.

\begin{equation*}
T_{R\circ R}^{j}(y,x)=\underset{z}{\vee }\left\{
T_{R}^{j}(y,z)\wedge T_{R}^{j}(z,x)\right\} \geq \underset{z}{\vee
}\left\{ T_{R\circ R}^{j}(y,z)\wedge T_{R\circ R}^{j}(z,x)\right\}
=T_{R^{2}\circ R^{2}}^{j}(y,x),
\end{equation*}
\begin{equation*}
I_{R\circ R}^{j}(y,x)=\underset{z}{\wedge }\left\{
I_{R}^{j}(y,z)\vee I_{R}^{j}(z,x)\right\} \leq \underset{z}{\wedge
}\left\{ I_{R\circ R}^{j}(y,z)\vee I_{R\circ R}^{j}(z,x)\right\}
=I_{R^{2}\circ R^{2}}^{j}(y,x)
\end{equation*}
and
\begin{equation*}
F_{R\circ R}^{j}(y,x)=\underset{z}{\wedge }\left\{ F(y,z)\vee
F_{R}^{j}(z,x)\right\} \leq \underset{z}{\wedge }\left\{ I_{R\circ
R}^{j}(y,z)\vee F_{R\circ R}^{j}(z,x)\right\} =F_{R^{2}\circ
R^{2}}^{j}(y,x)
\end{equation*}
Finally, the proof is valid.

\section{Acknowledgements}

The authors would like to thank the anonymous reviewer for their
careful reading of this research paper and for their helpful
comments.

\section{Conclusion}

In this paper, we have firstly defined the neutrosophic multi
relations(NMR). The NMR are the extension of neutrosophic soft
relation(NR)\cite{del-14-1} and intuitionistic multi relation
\cite{raj-13b}. Then, some notions such as; inverse, symmetry,
reflexivity and transitivity on neutrosophic multi relations are
studied. The future work will cover the application of the MNR in
decision making, pattern recogntion and in medical diagnosis.

\end{document}